\documentclass[reqno,11pt]{amsart}

\usepackage{color}

\newtheorem{theorem}{Theorem}[section]

 \theoremstyle{remark}

\theoremstyle{definition}

\newcommand\bR{\mathbb{R}}

\newcommand{\mysection}[1]{\section{#1}
      \setcounter{equation}{0}}

\begin{document}

\title[On a result of Jessica Lin]
{On a result of Jessica Lin}

\author{N.V. Krylov}
\thanks{The  author was partially supported by
 NSF Grant DMS-1160569}
\email{krylov@math.umn.edu}
\address{127 Vincent Hall, University of Minnesota,
 Minneapolis, MN, 55455}



 \begin{abstract}
In Theorem 1.5 of \cite{JL} an estimate from below
is given for nonnegative supersolutions of parabolic equations
in terms of the measure of the set where the right-hand side
 is less than $-1$. The given
proof contains many interesting details which are certainly
useful in many other situations. However, as long as
the statement of Theorem 1.5 of \cite{JL} is concerned,
it can be proved in a much shorter way, which we present here.
\end{abstract}

\maketitle

\mysection{Introduction}

The author of \cite{JL} presents a 
parabolic version of the lower bound for solutions
 established by Caffarelli, Souganidis, and Wang 
in \cite{4}, which was used in the error
estimates for stochastic homogenization of 
uniformly elliptic equations in random media \cite{3}. 
Although the general approach in \cite{JL} 
follows \cite{4}, it was necessary to develop a 
number of new arguments to handle the parabolic 
structure of the problem. One of the main technical tools
in \cite{JL} is Theorem 1.5, proved there by invoking  many 
facts which are very useful in many situations.
However, as long as
the {\em statement\/} of Theorem 1.5 of \cite{JL} is concerned,
it can be proved in a much shorter way, which we present here
expanding a too short to understand
  argument in Remark 3.4 of \cite{JL}
attributed to the author.
We  first consider  the case of linear equations and then
the general case.  

\mysection{Linear case}
                                   \label{section 10.26.1}

Let $\bR^{d}$ be a Euclidean space of points $x=(x^{1},
...,x^{d})$, $a(t,x)=(a^{ij}(t,x))$ a $d\times d$ matrix-valued
  function defined on $\bR^{d+1}=\bR\times\bR^{d}$.
Take a constant $\delta\in(0,1)$ and assume that for all
values of arguments and $\xi\in\bR^{d}$
$$
\delta|\xi|^{2}\leq a^{ij}\xi^{i}\xi^{j}\leq\delta^{-1}|\xi|^{2}.
$$

Introduce
$$
D_{i}=\frac{\partial}{\partial x^{i}},\quad D_{ij}=D_{i}D_{j},
\quad \partial_{t}=\frac{\partial}{\partial x^{i}},\quad
L=a^{ij}D_{ij}-\partial_{t}.
$$
Also set $B_{\alpha}=\{x\in\bR^{d}:|x|<\alpha\}$,
 $Q_{\alpha,\beta}= (0,\beta]\times B_{\alpha}$ and take a function
$u\in W^{1,2}_{d+1,loc}(Q_{1,1})\cap C(\bar{Q}_{1,1})$, which
is nonnegative on the parabolic boundary of $Q_{1,1}$.
Finally,  for a Borel subset  $\Gamma$ of $\bR^{d+1}$
denote by $|\Gamma|$ its Lebesgue measure.

\begin{theorem}
                                          \label{theorem 10.26.1}

Assume that $|\{Lu\leq -1\}\cap Q_{1,1}|\geq m|Q_{1,1}|$,
where $m\in(0,1)$. Then there exists a constant $N,\rho,\beta
\in(0,\infty)$ depending only on $\delta$ and $d$
and for any $\kappa\in(0,1)$
there exists a constant  $c  \in(0,\infty)$
depending only on $\delta,d$, and $\kappa$ such that
for $|x|\leq 1-\sqrt{1- \kappa},1\geq t\geq 1-m \kappa$ we have
\begin{equation}
                                              \label{10.26.2}
u(t,x)\geq cm^{\rho}
e^{-\beta/(\gamma m)}-N\|(Lu)_{+}\|_{L_{d+1}(Q_{1,1})},
\end{equation}
where $\gamma=1-\kappa$.
\end{theorem}

Proof. Observe that if the result
is true for $\kappa=1/2$, then it is obviously also
true for $\kappa\in(1,1/2)$
and one can even keep $c$ the same as for $\kappa=1/2$.
Therefore we concentrate on
$$
\kappa\in[1/2,1).
$$
Then, one can slightly shrink $Q_{1,1}$
and add to $u$ a small constant in order to preserve
the positivity of the new
function on the parabolic boundary of the new domain. After that passing to the limit 
in the corresponding version of \eqref{10.26.2}
 shows that without losing generality
we may and will assume that
$u\in W^{1,2}_{d+1}(Q_{1,1})\cap C(\bar{Q}_{1,1})$.
This fact,
for obvious reasons, allows us to   also assume that
$u$ and the coefficients of $L$ are infinitely differentiable
in $\bR^{d+1}$.
In that case, by adding to $u$ the solution
of the equation $Lv=-(Lu)_{+}$ in $Q_{1,1}$ with  
boundary data $v=-u$ ($\leq0$) on the parabolic boundary
of $Q_{1,1}$ and using the parabolic Alexandrov estimate,
we reduce the situation to the one in which $Lu\leq0$
in $Q_{1,1}$ and $u=0$ on the parabolic boundary
of $Q_{1,1}$. The last operation may yield a new function $u$
which is no longer infinitely differentiable in $
\bar{Q}_{1,1}$, but yet it will be continuous in $\bar{Q}_{1}$
and belong to $W^{1,2}_{p}(Q_{1,1})$ for any $p<\infty$.

After that introduce $\Gamma=
\{(t,x)\in Q_{1,1}:Lu(t,x)\leq-1\}$ and
denote by $w$ the solution of class $W^{1,2}_{d+1}
(Q_{1,2} )\cap C(\bar{Q}_{1,2} )$ of the problem
$$
Lw(t,x)=-I_{\Gamma}(t,x)I_{B_{1}}(x)I_{(0,T)}(t)
$$
with zero data on the parabolic boundary of $Q_{1,2} $,
where
$$
T=(1-m)(1+\gamma^{2}m) .
$$
As is easy to see 
$$
|\{(t,x)\in Q_{1,1} :
Lw(t,x)\leq -1\}|\geq m(1-m)|Q_{1,1}| \gamma^{2} 
\geq m |Q_{1,1}| \gamma^{2}/2.
$$
It follows by Theorem 4.1 of \cite{Kr_10} that
$w(2,0)\geq c_{0} m^{\rho}$, where $c_{0}=c_{0}(\delta,d,\kappa)
>0$, $\rho=\rho(\delta,d)<\infty$.

Furthermore, observe that if a point $t\geq1-\kappa m$,
then its distance to $(0,T)$ is greater than 
(recall that $\gamma\leq1/2$)
$$
1-\kappa m-T=\gamma m[1-\gamma (1-m) ]\geq \gamma m/2.
$$
Also if $|x|<1-\sqrt{\gamma}$, then the distance of $x$
to $\partial B_{1}$ is greater that $\sqrt{\gamma}\geq 
\sqrt{\gamma m/2}$.
By combining this with the fact that for $t\geq
T$ the function $w$ satisfies a homogeneous equation
by Theorem 3.6 of \cite{FSY}, for $|x|<\kappa$ and $t\geq  
1-\kappa m$, we get that
$$
w(2,0)\leq c(m\gamma/2)w(x,t),
$$
where $c(\mu)$ is a certain function
which depends only on $\mu,\delta $,
and $d$. Since by the maximum principle $w 
\leq u $ in $Q_{1,1}$ we conclude
$$
c(m\gamma/2)u(x,t)\geq c_{0} m^{\rho}.
$$

It only remains to figure out how $c(\mu)$ depends on $\mu>0$.
It follows from the three-line proof of Theorem 3.6 in \cite{FSY}
 that $c(\mu)$ can be taken as a constant
depending only on $\delta$ 
and $d$ times the Harnack constant, namely, the constant 
$d(\mu)$ which is such that
for any $(x,t), (y,s)\in Q_{1,1}$ for which $t\geq s+\mu\geq2\mu$,
 and $|x|,|y|\leq 1-\mu^{1/2}\nu$,
where $\nu=\nu(\delta,d)>0$, we have
\begin{equation}
                                         \label{8.23.1}
w(y,s)\leq d(\mu) w(x,t)
\end{equation}
for any $w\geq0$ satisfying the homogeneous equation in $Q_{1,1}$.

It turns out that it suffices to find $d(\mu)$ such that
\eqref{8.23.1} holds as long as $t\geq s+\mu\geq 2\mu $
 and $|x|,|y|\leq 1-\mu^{1/2} $. Indeed, if $\nu\geq1$,
then this is obvious. However if $\nu<1$, then by replacing
the condition $t\geq s+\mu\geq2\mu$,
 with $t\geq s+\mu\nu^{1/2}\geq2\mu\nu^{1/2}$,
  we again would get a stronger
result.

Thus, take $(x,t), (y,s)\in Q_{1,1}$
such that
$t\geq s+\mu\geq 2\mu $
 and $|x|,|y|\leq 1-\mu^{1/2} $,
 consider the straight segment joining them,
and define $n$ as the smallest integer such that
$|x-y|/n \leq\sqrt\mu/\sqrt{n}$. Obviously, $n\leq 5/\mu $
for small $\mu$.
Then we split the segment into $n$ parts of equal length, so that
for any two neighboring points $(x_{i},t_{i})$ and
$(x_{i+1},t_{i+1})$ we have $t_{i+1}-t_{i}=(t-s)/n$ ($\geq
\mu /n$) and $|x_{i+1}-x_{i}|=|x-y|/n \leq\sqrt\mu/\sqrt{n}$.
Now each couple of points $(x_{i},t_{i})$ and
$(x_{i+1},t_{i+1})$ is in a standard-shape cylinder to which
Harnack's inequality is applicable, so that
$$
w(x_{i},t_{i})\leq d_{0}w(x_{i+1},t_{i+1}),
$$
where $d_{0}=d_{0}(\delta,d)$,
implying that \eqref{8.23.1} holds with $d(\mu)
=d_{0}^{n}\leq \exp(\beta/\mu )$.
This proves the theorem.

\mysection{Nonlinear case}

Let $A$ be a closed set of $d\times d$ symmetric matrices
with eigenvalues in $[\delta,\delta^{-1}]$.
Define
\begin{equation}
                                                  \label{10.27.2}
F[v](t,x)=\max_{a\in A}a^{ij}D_{ij}v(t,x)-\partial_{t}v(t,x).
\end{equation}
We are going to deal with $L_{p}$-viscosity
supersolutions the definition of which
and their numerous properties can be found in \cite{CKS00}.
Here is a somewhat more detailed version of Theorem 1.5
of \cite{JL}.
\begin{theorem}
                                          \label{theorem 10.27.1}
Let $f\in L_{d+2}(Q_{1,1})$ and $u\in C(\bar{Q}_{1,1})$
be such that $u,f\geq0$ and $u$ is an $L_{d+2}$-viscosity
supersolution of the equation $F[v]=-f$ in $Q_{1,1}$.
Assume that $|\{f\geq1\}\cap Q_{1,1}|\geq m|Q_{1,1}|$,
where $m\in(0,1)$. Then
for $\kappa\in(0,1)$, $|x|\leq \kappa,1\geq t\geq 1-m \kappa$ we have
\begin{equation}
                                                  \label{10.27.1}
u(t,x)\geq cm^{\rho}
e^{-\beta/(\gamma m)},
\end{equation}
where $c,\rho,\beta,\gamma$ are taken from Theorem \ref{10.26.2}.
\end{theorem}

Proof. By Theorem 8.4 of \cite{CKS00}
the equation $F[v]=-f$ in $Q_{1,1}$ with zero condition
on the parabolic boundary of $Q_{1,1}$ has a unique
solution $v\in W^{1,2}_{d+2,loc}(Q_{1,1})\cap C(\bar{Q}_{1,1})$.
It also follows from \cite{CKS00} that $v\leq u$ in $\bar{Q}_{1,1}$.
Furthermore, observe that since
 $A$ is closed, the maximum in \eqref{10.27.2} is attained
for each $(t,x)$,
and there exists an operator $L$ as in Section \ref{section 10.26.1}
such that $Lv=-f$ in $Q_{1,1}$. After that it only remains to apply Theorem
\ref{theorem 10.26.1}. The theorem is proved.

{\bf Acknowledgment}. The author is sincerely
grateful to Hongjie Dong who found some
inconsistencies in the first draft of the paper.

\end{document}